\documentclass[11pt]{article}
\usepackage{amssymb, amsmath}
\usepackage{color, pdfcolmk}
\usepackage{graphics}
\usepackage{graphicx}
\usepackage{amssymb}
\usepackage{epsf}
\usepackage{verbatim}

\newtheorem{theorem}{Theorem}
\newtheorem{proposition}[theorem]{Proposition}
\newtheorem{lemma}[theorem]{Lemma}
\newtheorem{corollary}[theorem]{Corollary}

\newtheorem{definition}{Definition}
\newtheorem{example}{Example}

\newcommand \bB{\mathbb{B}}
\newcommand \bD{\mathbb{D}}

\newcommand \bR{{\mathbb{R}}}

\newcommand \eps{\varepsilon}
\newcommand \vs{\vec{s}}
\newcommand \vsst{\vec{s}^{\, *}}
\newcommand \vx{\vec{x}}
\newcommand \vy{\vec{y}}
\newcommand \vz{\vec{z}}
\newcommand \vg{\vec{\gamma}}

\newcommand \xss{x^{**}}
\newcommand \xs{x^*}
\newcommand \xm{x^\circ}
\newcommand \cQ{\mathcal{Q}}

\begin{document}

\title{Two classes of ODE models with switch-like behavior}

\author{Winfried Just\footnote{Corresponding author. 
Department of Mathematics, Ohio University, Athens, OH 45701, U.S.A.  E-mail: mathjust@gmail.com}, Mason Korb, Ben Elbert, Todd Young\footnote{T.Y.\ and this work were partially supported by the NIH-NIGMS grant R01GM090207.}\\
Ohio University}

\maketitle

\begin{abstract}
In cases where the same real-world system can be modeled both by an ODE system~$\bD$ and a Boolean system~$\bB$ it is of interest to identify conditions under which the two systems will be consistent, that is, will make qualitatively equivalent predictions. In this note we introduce two broad classes of relatively simple models that provide a convenient framework for studying such questions.  In contrast to the widely known class of Glass networks, the right-hand sides of our ODEs are Lipschitz-continuous.  We prove that under suitable assumptions about~$\bB$ consistency between $\bD$ and $\bB$ will be implied by sufficient separation of time scales in one class of our models while it may fail in the other class. These results appear to point to more general structure properties that favor consistency between ODE and Boolean models.
\end{abstract}

\section{Introduction}

The dynamics of many real-world, in particular biological, systems exhibits switchlike behavior. Examples include gene regulatory networks, where at any given time a gene may be expressed (switched on) or not expressed (switched off) and neuronal networks, where  at any given time a neuron may or may not fire.  Mathematical models of such systems can often take the form of a differentiable flow~$\bD$ or a Boolean dynamical system~$\bB$.  The former type of model usually incorporates more relevant details and appears biologically more realistic, while the latter type may be easier to analyze, at least by simulations. The literature describes a vast number of cases where Boolean and other discrete models have been successfully used to make realistic predictions about real-world systems; for example, several chapters of~\cite{Raina} review such success stories in biology.

If the same real-world system can be modeled both as a differentiable flow~$\bD$ and a Boolean system~$\bB$,
the question naturally arises how the two models are related to each other.
This question has several distinct aspects, two of which will be discussed in this note: translation between the two types of systems and consistency.

\subsection{Translation between flows and Boolean systems}\label{subsec:translation}

Differentiable flows~$\bD$ are induced by ODEs

\begin{equation}\label{eqn:ODE-general}
\dot{\vx} = g(\vx).
\end{equation}

While formally different mathematical objects, for convenience we will sometimes refer to  the flow~$\bD$ and the ODE~(\ref{eqn:ODE-general}) as if they were synonymous.  The state space of a differentiable flow is  a subset~$M$  of some~$\bR^N$, typically a compact manifold.

The state space of an $n$-dimensional Boolean system is the set~$\{0,1\}^n$ of $n$-dimensional Boolean vectors~$\vs$.
The dynamics is given by a Boolean function $f: \{0,1\}^n \rightarrow \{0,1\}^n$ so that

\begin{equation}\label{eqn:Booldyn}
\vs(\tau+1) = f(\vs(\tau)).
\end{equation}

While a Boolean system is formally a pair~$\bB = (\{0,1\}^n, f)$, the state space is already implied by the updating function~$f$ as its domain, so we will sometimes use~$f$ itself as shorthand for specifying~$\bB$.

Translating a flow~$\bD$ into a Boolean system~$\bB$ requires at the minimum construction of a map
$S: M \rightarrow \{0,1\}^n$ that we will call the \emph{discretization} of~$M$.  Typically, this map is constructed by discretizing the values $x_i$ for~$n$ among the variables in~$M$ and assigning to them the Boolean value~$s_i = 0$ if $x_i$ is below a threshold~$\Theta_i$ and $s_i = 1$ if $x_i$ is above this threshold.  In practical applications it may be far from obvious how to choose the thresholds~$\Theta_i$, see~\cite{Discretization}, but this issue is beyond the scope of this note. We want to point out though that not all variables in~$M$ need to be discretized. For example, in models of gene regulation one may choose to incorporate only discretized values of mRNA concentrations but not of their translated proteins into the Boolean model.  Similarly, the Boolean models constructed in~\cite{TAWJ} have only one Boolean variable per neuron,   while the state of each each neuron in the ODE model is characterized both by the cross-membrane voltage and a second, so-called gating variable.

Translating the ODEs~(\ref{eqn:ODE-general}) into a Boolean updating function~$f$ is less straightforward than discretizing some variables, except in special cases when we have \emph{strong consistency} (see Section~\ref{sec:consistency}), or when $f$ itself is inherent in the construction of~$\bD$. In practice, the updating functions of Boolean models are usually based directly on the available empirical knowledge of the real-world system. Fortunately, for the results reported here, it does not matter how~$f$ is constructed. Therefore we will here somewhat informally assume that a given flow $\bD$ has been translated into a Boolean system~$\bB$  by means of a discretization~$S$ and some unspecified ``natural'' choice of~$f$ and call the resulting $\bB$ a \emph{(natural) Boolean approximation of~$\bD$.} We need to emphasize though that our usage of these phrases does not imply any consistency whatsoever between a given flow and its natural Boolean approximation.

Wile at least some differentiable flows~$\bD$ have unique natural Boolean approximations, translation from Boolean into ODE systems cannot be a one-to-one procedure.  However, for some types of real-world systems we may have a fairly good idea about a general class of ODE models that give a reasonably good description of the underlying mechanisms. Being able to identify a subset or a single representative of this class that will have a given Boolean system as its approximation can be a powerful tool for model selection in cases where the available data allow for construction of Boolean models.  The caveats about non-uniqueness nonwithstanding, we like to think about this process as translating Boolean into ODE models.

For example, \emph{Glass networks} are an extensively studied class of simple ODE models for gene regulatory networks~\cite{Edwards}, \cite{Glass1}--\cite{Glass4}.
 In these networks, the concentration $x_i$ of the $i$-th mRNA (or gene product, if one prefers to think about them this way) is modeled by a DE

\begin{equation}\label{eqn:GlassODEs}
\dot{x}_i = F_i(\vx) - \lambda_i x_i,
\end{equation}
where $F_i$ is a piecewise constant (usually discontinuous)  production term and  $-\lambda_i x_i$ is the decay term.  Under certain technical assumptions  a Glass network has a natural Boolean approximation.

Similarly, Terman \emph{et al.}  \cite{TAWJ} constructed a class of excitatory-inhibitory neuronal networks whose
ODE models have natural Boolean approximations and showed that every Boolean system can be translated into a model in this class.  Together with the theorem about consistency that was proved in~\cite{TAWJ}, this makes both the continuous and discrete dynamical systems promising candidates for modeling actual neuronal networks that exhibit a phenomenon called \emph{dynamic clustering.} As in Glass networks, some ODEs of~\cite{TAWJ} have right-hand sides with discontinuities.

Here we introduce two other classes of differentiable flows with natural Boolean approximations so that every Boolean system can be translated into a set of flows in each class. These classes are very broad, relatively easy to study, and the ODEs in our classes have Lipschitz-continuous right-hand sides. While these classes are not directly related to any applications that we are aware of, we believe they provide a convenient set of toy models for exploring the relation between differentiable flows and their Boolean approximations, such as conditions for consistency.   In particular, the study of these classes may allow us to identify which type results about Glass networks or from~\cite{TAWJ} do require discontinuities in the right-hand sides of the ODEs, such as in~(\ref{eqn:GlassODEs}), and which ones do not.

\subsection{Consistency}\label{subsec:continuity}

If a flow $\bD$ and a Boolean system $\bB$ are to function as useful models of the same real-world system, they should be \emph{consistent} in the sense that they make equivalent predictions.  In particular, we would like to assure that for all trajectories of~$\bD$ that start in a sufficiently large subset~$U$ of~$M$ the next Boolean state of the system will be what $f$, or at least one of the coordinates of~$f$, predicts for the current Boolean state $S(\vx(t))$. While quite intuitive, this notion is not easy to formalize; in fact, it admits a number of nonequivalent formal definitions that seem plausible. The choice of the most appropriate one may depend on the context (see Section~\ref{sec:discussion} for a brief discussion of one example). The problem is that while each state $\vs(\tau)$ of~$\bB$ has a well-defined successor state $\vs(\tau+1) = f(\vs(\tau))$, time in the flow~$\bD$ increases continuously and one has to pick a time $t^+ > t$ for which $S(\vx(t^+))$ will be considered the ``next'' Boolean state after time~$t$.
Here we will study two such formal and quite natural notions; we call them \emph{strong consistency} and \emph{consistency.}

\subsection{Outline of the paper and main results}\label{subsec:outline}

In Section~\ref{sec:translation} we will describe methods for translating, in the sense of Subsection~\ref{subsec:translation}, a Boolean system $\bB = (\{0,1\}^n, f)$ into differentiable flows $\bD_1(f, \vg)$ and $\bD_2(f, \vg)$ with state spaces~$I^n$ and~$I^{2n}$ respectively, where $I$ will be a certain compact interval. Each real variable in systems~$\bD_1(f, \vg)$ will correspond to a separate Boolean variable; while for systems~$\bD_2(f, \vg)$ only the variables~$x_i$ for $i \in \{1, \ldots, n\}$ will have Boolean discretizations.  We call these variables the \emph{signature variables} as they determine the Boolean state of the flow.  Each signature variable~$x_i$ receives input only from $x_{i+n}$ whereas $x_{i+n}$ receives input only from the signature variables $x_1, \ldots , x_n$.

One can think about the systems~$\bD_2(f, \vg)$ as describing the dynamics of \emph{agents} $X_i = \{x_i, x_{i+1}\}$ whose internal states are represented by the $x_i$s and who receive signals about the state of the system as a whole through $x_{i+n}$.  Thus we will refer to the variables~$x_{i+n}$ as the \emph{signaling variables.}  The parameter~$\vg$ will allow us to control the degree of separation of time scales between variables of the flows; in particular, between the dynamics of the signature variables and the signaling variables.  Let us mention that
$\bD_1(f, \vg)$ and $\bD_2(f, \vg)$ also depend on another parameter~$\cQ$ that can be thought of as the particular translation method or \emph{conversion scheme.}  Most of our results allow a lot of flexibility in the choice of~$\cQ$ but assume~$\cQ$ fixed; therefore we will often suppress this parameter in our notation.

In Section~\ref{sec:consistency} we formally define the notions of strong consistency and consistency. In Section~\ref{sec:Boolean} we define the two important classes of one-stepping and monotone-stepping Boolean systems and show that strong consistency is possible only if the Boolean system is one-stepping.
In Section~\ref{sec:examples} we  give some examples of (strong) consistency or lack thereof for very simple Boolean systems.
In Section~\ref{sec:Consistency-Stepping} we prove the two main results of this paper, namely that sufficient separation of time scales guarantees strong consistency of~$\bD_2(f, \vg)$ with~$f$ for one-stepping~$f$  and consistency of~$\bD_2(f, \vg)$ with~$f$ for monotone-stepping~$f$.
In Section~\ref{sec:discussion} we discuss the relation of our results to similar work in the literature and we outline some directions of future research.

\subsection{Notation}\label{subsec:notation}

Our mathematical notation is mostly standard. The set $\{1, \ldots, n\}$ will be denoted by~$[n]$. The cardinality of a set~$X$ will be denoted by~$|X|$.
Time in flows will be denoted by~$t$. Since most Boolean systems are not time-reversible, we will only be interested in forward trajectories of the corresponding flows and will always assume that $t \in [0, \infty)$.  States in flows will be denoted by $\vx$, with $\vx(t)$ denoting the state at time~$t$ of the trajectory with initial state~$\vx(0)$. Time in Boolean systems will be denoted by~$\tau$ and will be assumed to be a nonnegative integer.
The components of a real-valued vector~$\vx$ will be denoted by~$x_i$; similarly, the components of a Boolean vector~$\vs$ will be denoted by~$s_i$. Trajectories in a Boolean system~$\bB = (\{0,1\}^n, f)$ are sequences  $\{\vs(\tau)\}_{\tau = 0}^\infty$ such that $\vs(\tau) = f^\tau(\vs(0))$ for all~$\tau$. The \emph{Hamming distance}
$H(\vs, \vs^{\, *})$ between two Boolean vectors is the number of indices~$i$ such that $s_i \neq s^*_i$.

\section{The flows $\bD_1(f, \vg)$ and $\bD_2(f, \vg)$}\label{sec:translation}

\subsection{General definition}\label{subsec:generalD12}

Fix~$n$ and let $s \in \{0,1\}^n$.  Let $N \geq n$ and define sets $W^s_i \subset \bR^N$ for $i \in [n]$ as follows:

\begin{equation}\label{eqn:Wsi-def}
W^s_i = \{\vx: \ x_i < -1\} \ \mbox{if} \ s_i = 0, \qquad
W^s_i = \{\vx: \ x_i > 1\} \ \mbox{if} \ s_i = 1.
\end{equation}

Let $W^s = \bigcap_{i \in [n]} W^s_i$.

A \emph{continuous conversion} of a Boolean updating function $f: \{0,1\}^n \rightarrow \{0, 1\}^n$ is a continuous function
$$Q = (Q_1, \dots , Q_n): \bR^N \rightarrow [0,1]^n$$
such that

\begin{equation}\label{eqn:Q-reflects-Boolean}
\forall s \in \{0, 1\}^n \, \forall \vx \in W^s \, \forall i \in [n] \ Q_i(\vx) = f_i(\vs),
\end{equation}
where $f_i$ is the $i$-th component of~$f=(f_1, \dots , f_n)$.
For example, the piecewise linear function

\begin{equation}\label{eqn:L}
L(x)=
\begin{cases}
 0 	 	&		    \mbox{if }  x \leq -1, \\
 .5(x+1) &   \mbox{if } -1 < x < 1,  \\
 1       &     \mbox{if } x \geq 1,
\end{cases}
\end{equation}
allows us to construct continuous conversions $L(x_i)$ for the projections $\pi_i(s) = s_i$ and $L^n := (L(x_1), \dots , L(x_n))$ of the identity function  $f_{id}(s) = s$.

A \emph{conversion scheme} is a map~$\cQ$ that assigns to each Boolean function a continuous conversion~$Q(f)$.

Now let $\cQ$ be a conversion scheme, let $f: \{0,1\}^n \rightarrow \{0, 1\}^n$ be a Boolean function, and let $\vg$ be an $n$-dimensional vector of positive reals. Define an ODE systems $\bD_1(f, \vg, \cQ)$ (or simply $\bD_1(f, \vg)$ if $\cQ$ is assumed fixed or implied by the context) by

\begin{equation}\label{eqn:standardform}
\dot{x}_i=\gamma_i (g(x_i)+ 6 Q_i(\vx)),
\end{equation}
where

\begin{equation}\label{eqn:gi}
g(x)=3 x- x^3-3,
\end{equation}
and $Q_i$ is the $i$-th coordinate of $Q = Q(f)$.

While we assume here for simplicity  that the parameters $\gamma_i$ are constants, we want to remark that the arguments in this paper can be generalized to the case when these values are allowed to
depend on the state or even change over time in a nonautonoumous way, as long as they are all bounded and bounded away from zero, that is, if there are constants $M > m > 0$ such that $m < \gamma_i(\vx, t) < M$ for all~$i, \vx$.

In order to see how the systems $\bD_1(f, \vg)$ behave,
consider the family of flows generated by the one-dimensional ODEs

\begin{equation}\label{eqn:ODE-g-h}
\dot{x} = \gamma (g(x)+ h),
\end{equation}
where we consider
both $\gamma$ and $h$ as potential bifurcation parameters.

The only bifurcations with respect to~$\gamma$ occur at the bifurcation value~$\gamma^* = 0$. These will be of no interest to us, since we will assume that $\gamma > 0$.  Under this assumption $\gamma$ does not alter the qualitative behavior of~(\ref{eqn:ODE-g-h}) but controls the speed at which trajectories are being traversed, with small $\gamma$ corresponding to slow change. Thus suitable choices of~$\vg$ will allow us to study the effect of separation of time scales.

Next note that $g$ has two local extrema, a local minimum~$g(-1) = - 5$  and a local maximum~$g(1) = -1$.  Thus for $\gamma > 0$ the family of flows defined by~(\ref{eqn:ODE-g-h}) undergoes two saddle-node bifurcations with respect to~$h$:
For $h < h^{*} = 1$, there will be exactly one globally stable equilibrium $\xs < -1$; for $1 = h^{*} < h < h^{**} = 5$ there will be two locally stable equilibria~$\xs < -1$ and $\xss > 1$ as well as one unstable equilibrium~$\xm \in (0,1)$; and
for $h > h^{**} = 5$, there will be exactly one globally stable equilibrium $\xss > 1$.
This translates into the context of~(\ref{eqn:standardform}) as follows.

\begin{proposition}\label{prop:bifurcations}
Assume~$Q_i$ is a constant.  If  $Q_i < 1/6$, then~(\ref{eqn:standardform}) will have exactly one globally asymptotically stable equilibrium $\xs < -1$; if $1/6 < Q_i < 5/6$ there will be two locally stable equilibria~$\xs < -1$ and $\xss > 1$ as well as one unstable equilibrium~$\xm \in (0,1)$; and if $Q_i > 5/6$, there will be exactly one globally asymptotically stable equilibrium $\xss > 1$.
\end{proposition}

Of course in our actual systems  $\bD_1(f, \vg)$ and $\bD_2(f, \vg)$ the values of $Q_i$ will change along trajectories; much of our work in this note will be concerned with studying  the consequences of such changes.

For $\cQ, f$ as above and a $2n$-dimensional vector $\vg$ of positive reals, we define an ODE systems $\bD_2(f, \vg, \cQ)$ (or simply $\bD_2(f, \vg)$ if $\cQ$ is assumed fixed or implied by the context) with variables $x_i, x_{i+n}$ for $i \in [n]$ by

\begin{equation}\label{eqn:standardform2}
\begin{split}
\dot{x}_i &=\gamma_i (g(x_i)+ 6 L(x_{i+n})),\\
\dot{x}_{i+n} &=\gamma_{i+n} (g(x_{i+n})+ 6Q_{i}(x_1, \dots , x_n)),
\end{split}
\end{equation}
where $Q_i$ is the $i$-th coordinate of $Q = Q(f)$ and $L, g$ are as in~(\ref{eqn:L}), (\ref{eqn:gi}).

In a sense, equations~(\ref{eqn:standardform2}) are special cases of equations~(\ref{eqn:standardform}). To see this, let  $f=(f_1, f_2, \dots, f_n): \{0, 1\}^n \rightarrow \{0,1\}^n$ be given. For each $i\in [n]$ we define an auxiliary function $c_{i}(\vec{s})= s_{n+i}$ that copies the value of variable number~$n+i$ to variable number~$i$.
Extend~$f$ to a Boolean function $f^+: \{0, 1\}^{2n} \rightarrow \{0, 1\}^{2n}$ given by

\begin{equation}\label{eqn:f+def}
f^+= (c,f) = (c_1 , \dots, c_n, f_1,\dots, f_{n}).
\end{equation}
Then assuming a conversion scheme~$\cQ$ with $Q(\pi_{i+n}) = L(x_{i+n})$ the system~(\ref{eqn:standardform2}) is the same as the system~(\ref{eqn:standardform}) that defines  $\bD_1(f^+, \vg, \cQ)$.

Since  $L$ and $Q_{i}$ take only values in the interval $[0,1]$, we make the following observation (see~\cite{LowDimEx16} for a proof).

\begin{proposition}\label{prop:statespace}
Let $f$ be an $n$-dimensional Boolean function and let $\vg$ denote any vector of positive reals of suitable dimension.  Let $x^- \approx -2.1038$ be the unique root of the polynomial $g(x) = 3x - x^3 - 3$ and let
$x^+ \approx 2.1038$ be the unique root of the polynomial $g(x) + 6 = 3x - x^3 + 3$. Then $[x^-, x^+]^n$ is a forward-invariant set in $\bD_1(f, \vg)$ and $[x^-, x^+]^{2n}$ is a forward-invariant set in $\bD_2(f, \vg)$.
\end{proposition}

Thus for our purposes we will consider $M = [x^-, x^+]^n$ as the state space of the flow $\bD_1(f, \vg)$ and $M = [x^-, x^+]^{2n}$ as the state space of the flow $\bD_2(f, \vg)$. We define discretizations of individual variables in these flows by

\begin{equation}\label{eqn:S-i}
S_i(\vx)=
\begin{cases}
 0 	 			&  \text{if $x_i \leq 0$, } \\
 1        &  \text{if $x_i > 0$,}
\end{cases}
\end{equation}
and define $S: M \rightarrow \{0,1\}^n$ as

\begin{equation}\label{eqn:S}
S(\vx)= \prod_{i \in [n]} S_i(\vx).
\end{equation}
In particular, $S$ maps each of the sets $W^s$ to~$s$.

Note that even though the systems~$\bD_2(f, \vg)$ have dimension~$2n$,  both in the case of~$\bD_1(f, \vg)$ and of~$\bD_2(f, \vg)$ the Boolean state $S(\vx)$ depends only on the variables $x_1, \ldots , x_n$.  Let us call these variables the \emph{signature variables} of the ODE system. According to~(\ref{eqn:standardform2}), in the systems~$\bD_2(f, \vg)$ none of these signature variables is directly influenced by signature variables; all interactions  are mediated by the variables $x_{n+1}, \dots , x_{2n}$.    Therefore we will refer here to the variables $x_{n+1}, \dots , x_{2n}$ as the \emph{signaling variables.}

\subsection{Examples of conversion schemes}\label{subsec:convert-schemes}

Here we give examples of several specific conversion schemes~$\cQ$.  All of them represent the functions~$Q_i$ as compositions

\begin{equation}\label{eqn:Q-P-L}
Q_i(\vx) = P_i(L(x_1), \dots , L(x_n)),
\end{equation}
where $P_i: [0,1]^n \rightarrow [0,1]$ is a continuous function. Such
$P = (P_1, \ldots , P_n)$
will result in a continuous conversion~$Q$ of a Boolean function $f: \{0, 1\}^n \rightarrow \{0, 1\}^n$
that satisfies~(\ref{eqn:Q-reflects-Boolean}) as long as

\begin{equation}\label{eqn:P-reflects-Boolean}
\forall s \in \{0, 1\}^n \,  \forall i \in [n] \ P_i(\vs) = f_i(\vs).
\end{equation}

Wittman \emph{et al.}~\cite{Wittman} give an algorithm for constructing  polynomial functions of minimal possible degree that satisfy~(\ref{eqn:P-reflects-Boolean}) and prove their uniqueness.  Together with~(\ref{eqn:Q-P-L}) these minimal-degree polynomials define a conversion scheme~$\cQ_W$.

While polynomial functions of minimal degree have desirable properties from our point of view (see
Subsection~\ref{subsec:constsec1}) and may be in some ways optimal, we want to allow also for other possible conversion procedures.

A drawback of~$\cQ_W$ is that the algorithm of~\cite{Wittman} requires evaluation of~$f_i(s)$ at $2^n$ data points, even if $f_i$ is given by a very simple Boolean formula. But for simple~$f_i$ it may be much easier to construct a suitable polynomial $P_i$.  For example, the following choices of~$P_i$ satisfy~(\ref{eqn:P-reflects-Boolean}):

\begin{itemize}
\item If $f_i = s_j \wedge s_k$, let  $P_i = x_jx_k$.
\item If $f_i = \neg s_j$, let  $P_i = 1 - x_j$.
\item If $f_i = s_j \vee s_k$, let  $P_i = x_j + x_k - x_jx_k$.
\item If $f_i = s_j \oplus s_k$, where $\oplus$ denotes exclusive or, let  $P_i = x_j + x_k - 2x_jx_k$.
\end{itemize}

If $f_i$ is given in conjunctive or disjunctive normal form, then one can use the first three observations above to recursively construct~$P_i$, which in conjunction with~(\ref{eqn:Q-P-L}) defines recursive conversion
schemes~$\cQ_{Rc}, \cQ_{Rd}$.  Similarly, since the operation~$\oplus$ is the same as addition in the field
$\mathbb{F}_2$, the first and fourth of the above observations can be used to define a recursive conversion
schemes~$\cQ_{RF}$ for Boolean functions that are given as polynomials over the field~$\mathbb{F}_2$ in the form

\begin{equation}\label{eqn:F2-poly}
f_i(s) = \prod_{j\in J_1} s_j \oplus \prod_{j\in J_2} s_j \oplus \dots \oplus \prod_{j\in J_k} s_j,
\end{equation}
where $J_1, J_2, \dots , J_k$ are subsets of~$[n]$.

In all three cases, if the relevant  expressions for $f_i$ contain few terms, these recursive conversion schemes allow for faster computation of  polynomials~$P_i$ than the algorithm of~\cite{Wittman} for~$\cQ_W$, but are not guaranteed to return polynomials of minimal degree.

Suppose $f_i$ is given as a  polynomial over~$\mathbb{F}_2$ as in~(\ref{eqn:F2-poly}).
  If $k$ is relatively large, then using the recursive conversion method~$\cQ_{RF}$ becomes unwieldy.  However, one can  find a quick and easy conversion as follows.  First replace~(\ref{eqn:F2-poly}) with a polynomial function $u: \bR^n \rightarrow \bR$ defined by

\begin{equation}\label{eqn:F2-to-real-poly}
u(\vx) = \prod_{j\in J_1} x_j + \prod_{j\in J_2} x_j + \dots + \prod_{j\in J_k} x_j.
\end{equation}

Then define

\begin{equation}\label{eqn:u-to-[01]}
P_i(\vx) =  0.5 - 0.5 \cos(\pi u(\vx)).
\end{equation}

The resulting function $P_i$ is no longer a polynomial, but it is continuous, even analytic, and maps $[0,1]^n$ into $[0,1]$.  According to~(\ref{eqn:F2-poly}), for $\vs \in \{0,1\}^n$ the value $f_i(\vs)$ is~0 if an even number of the products~$\prod_{j\in J_\ell} s_j$ evaluates to~1 and  is~1 if an odd number of the products~$\prod_{j\in J_\ell} s_j$ evaluates to~1.  In other words, $f_i(S(\vx))=0$ if $u(\vx)$ is an even integer and $f_i(S(\vx))=1$ if $u(\vx)$ is an odd integer. Now~(\ref{eqn:P-reflects-Boolean}) follows from~(\ref{eqn:u-to-[01]}).

This defines yet another conversion scheme~$\cQ_a$; we call it \emph{arithmetic conversion.}  It has the advantage of allowing very fast computation of the conversion for Boolean functions that are represented as polynomials over~$\mathbb{F}_2$. A potential disadvantage is that $P_i$ will in general take the values~$0$ and~$1$ at many points in the interior of~$[0,1]^n$, which may introduce additional equilibria of~$\bD_1(f, \vg, \cQ_a)$ or~$\bD_2(f, \vg, \cQ_a)$ that have no counterparts in terms of~$f$. While we have not studied in detail whether these equilibria may adversely affect consistency in some cases, the results of Section~\ref{sec:Consistency-Stepping} show that they do not destroy consistency in~$\bD_2(f, \vg, \cQ_a)$ under the assumption of sufficient separation of time scales.

Finally, we want to point out that for all conversion schemes described here the right-hand sides of the ODEs~(\ref{eqn:standardform}) and~(\ref{eqn:standardform2}) will be Lipschitz-continuous.

\subsection{Software for conversion and exploration}\label{subsec:software}

We wrote a software package \emph{Boolean-Continuous GUI}  that runs under MatLab  and is capable of constructing, for  user-specified~$f$ and~$\vg$, ODE systems~$\bD_1(f, \vg, \cQ)$ and $\bD_2(f, \vg, \cQ)$, where $\cQ \in \{\cQ_{RF}, \cQ_a\}$. The software also allows for simulating trajectories in these flows, tracking the corresponding Boolean states, and detecting certain additional features. Both the software and its user-manual are available from the authors upon request.

\section{Consistency and Strong Consistency: Definitions}\label{sec:consistency}

Suppose~$\bD$ is a flow on a forward-invariant  set $M \subseteq \bR^N$, let $\bB = (\{0,1\}^n, f)$ be an $n$-dimensional Boolean system, and let $S: M \rightarrow \{0, 1\}^n$ be a fixed discretization. Let~$S_i$ denote the $i$-th component of~$S$ and let

$$B = \bigcup \{bd(S_i^{-1}\{0\}\}.$$

We will say that the trajectory of~$\vx = \vx(0)$ has the \emph{transversality property} if there exists an increasing  sequence~$\{t_k\}$  of nonnegative reals,
where $0 \leq k < K \leq \infty$, such that

\begin{itemize}
\item[(a)] $t_0 = 0$.
\item[(b)] For all $t \in [t_0, t_1)$ we have $\vx(t) \notin B$.  In particular, $S(\vx(t))$ is constant for $t \in [t_0, t_1)$.
\item[(c)] For all $k + 1 < K$ and all $t \in (t_k, t_1)$ we have $\vx(t) \notin B$.  In particular, the function $S(\vx(t))$ is constant for $t \in (t_k, t_{k+1})$.
\item[(d)] For all $0 < k < K$ we have $\lim_{t \rightarrow t_k^-} S(\vx(t)) \neq \lim_{t \rightarrow t_k^+} S(\vx(t))$.
\item[(e)] If $K = \infty$, then $\lim_{k \rightarrow \infty} t_k = \infty$.
\end{itemize}

If $\vx$ does have the transversality property, then the sequence~$\{t_k\} = \{t_k(\vx) \}$ that satisfies~(a)--(e) above is uniquely defined.  We call it \emph{the sequence of switching times for~$\vx$.}  Given this sequence we can associate with~$\vx$ a sequence of Boolean states

\begin{equation}\label{eqn:Boolseq-vx}
\begin{split}
\vs^{\, \vx}(\tau) &= S\left(\vx\left(\frac{t_{\tau+1} - t_{\tau}}{2}\right)\right) \quad \mbox{if} \quad \tau+1 < K,\\
\vs^{\, \vx}(\tau) &= S(\vx(t_{K-1} + \tau + 1)) \quad \mbox{if} \quad \tau+1 \geq K.
\end{split}
\end{equation}

Recall the definition of Boolean dynamics~(\ref{eqn:Booldyn}).  This be written in the form of~$n$ updating functions for the components~$f_i$ of $f = (f_1, \ldots , f_n)$ as follows:

\begin{equation}\label{eqn:regs}
s_i(\tau+1) = f_i(s_1(\tau),\ldots,s_n(\tau)).
\end{equation}

A \emph{fixed point} of a Boolean system~$\bB$ (or of~$f$) is a state $s = (s_1, \ldots , s_n)$ such that $f_i(s) = s_i$ for all $i \in [n]$.

\begin{definition}\label{def:consistency} Let $\bD, M, f$ be as above.
\smallskip

\noindent
(i) We say that the trajectory of~$\vx \in M$ is \emph{strongly consistent} with $\bB$ if it has the transversality property and for all~$\tau$

\begin{equation}\label{eqn:strong-cons}
\vs^{\, \vx}(\tau+1) = f(\vs^{\, \vx}(\tau)).
\end{equation}

\smallskip

\noindent
(ii) We say that the trajectory of~$\vx \in M$ is \emph{consistent} with $\bB$ if it has the transversality property and for all~$\tau$

\begin{equation}\label{eqn:consis}
\begin{split}
\forall i\ ( s_i^{\vx}(\tau + 1) &=  f_i(s^{\vx}(\tau)) \ \vee \ s_i^{\vx}(\tau + 1) = s_i^{\vx}(\tau))\\
\vs^{\, \vx}(\tau + 1) & =  \vs^{\, \vx}(\tau) \ \mbox {iff} \ \vs^{\, \vx}(\tau) \ \mbox{is a fixed point of}\  \bB.
\end{split}
\end{equation}

\noindent
(iii) We say that the flow~$\bD$ is (strongly) consistent with~$\bB$ on $U \subseteq M$ if every trajectory that starts in~$U$ is (strongly) consistent with~$\bB$.

\smallskip

\noindent
(iv) We say that the flow is (strongly) consistent with~$\bB$  if it is (strongly) consistent with~$\bB$ on an open set~$U \subseteq M$ that is \emph{universal} in the sense that $S$ maps $U$ onto~$\{0,1\}^n$.
\end{definition}

Instead of the phrase ``consistent with~$\bB$'' we will often write ``consistent with~$f$'' for convenience.

\section{One-Stepping and Monotone-Stepping Boolean Systems}\label{sec:Boolean}

\begin{definition}\label{def:one-stepping}
Let $\bB = (\{0,1\}^n, f)$ be a Boolean system, and let $\vs\in \{0,1\}^n$.

\smallskip

\noindent
(i) We say that \emph{the trajectory $\{\vs(\tau)\}_{\tau = 0}^\infty$ of $\vs =\vs(0)$ is one-stepping} if
\begin{equation}\label{eqn:one-step-trajs}
\forall \tau \ H(\vs(\tau), \vs(\tau+1)) \leq 1.
\end{equation}

\noindent
(ii) We say that \emph{$\bB$ is one-stepping} if all its trajectories are one-stepping.
\end{definition}

\begin{example}\label{ex:fixed-vs-onestep}
(i) Every Boolean trajectory that starts at a fixed point is one-stepping.

\smallskip

\noindent
(ii) Let $I \subseteq [n]$ and define $f: \{0,1\}^n \rightarrow \{0,1\}^n$ by letting
$f_i(s) = s_{i-1}$ if $i \notin I$ and  $f_i(s) = 1- s_{i-1}$ if $i \in I$, where the subscript $i-1 =0$ gets treated as $i -1 = n$. If $|I|$ is odd, then $f$ defines a Boolean system~$\bB$ with exactly one one-stepping orbit that has length~$2n$.
For a proof see~\cite{Elbert} and pages 79--81 of~\cite{MasonThesis}.
\end{example}

\begin{example}\label{ex:one-stepping-system}
Let $\bB = (\{0,1\}^n, f)$ be any Boolean system.  Define a Boolean system  $\bB^\circ = (\{0,1\}^n, f^\circ)$ by letting $f^\circ_i(\vs) = f_i(\vs)$ if $f_j(\vs) = s_j$ for all $j$ with $1 \leq j < i$ and $f^\circ_i(\vs) = s_i$ otherwise. Then $\bB^\circ$ is one-stepping.
\end{example}

\begin{proposition}\label{prop:mono+cons-implies-strong}
Suppose $\bD$ is a flow on~$M$ and the trajectory of~$\vx \in M$ is consistent with a Boolean function~$f$ for a given discretization.  If the trajectory of $S(\vx)$ under iterations of~$f$ is one-stepping, then the trajectory of~$\vx$ is strongly consistent with~$f$.
\end{proposition}

\textbf{Proof:} By the assumption of consistency, the trajectory of~$\vx$ has the transversality property and the corresponding sequence of Boolean states
$\{\vs(\tau)\}_{\tau=0}^\infty = \{\vs^{\, \vx} (\tau)\}_{\tau=0}^\infty$ satisfies~(\ref{eqn:consis}).  By  assumption, the Boolean trajectory $\{f^\tau(\vs(0))\}_{\tau=0}^\infty$
 has the property that $f^\tau(\vs(0))$ differs from
$f^{\tau+1}(\vs(0))$ in at most one coordinate, and it follows from~(\ref{eqn:consis}) by induction over~$\tau$ that
$f^{\tau}(\vs(0)) = \vs(\tau)$, which in turn implies~(\ref{eqn:strong-cons}). $\Box$

\bigskip

While Proposition~\ref{prop:mono+cons-implies-strong} indicates that one-stepping trajectories in Boolean approximations of differentiable flows are promising candidates for strong consistency, the next result shows that these are in a sense the \emph{only} candidates for seeing strong consistency in a neighborhood.  This is a well-known result, but not in this terminology, so we include a proof for completeness. First we need some terminology and a lemma.

\begin{definition}\label{def:BoolODedef}
Let $M$ be a compact $m$-dimensional topological manifold with boundary for some~$n$, and let $S: M \rightarrow \{0,1\}^n$ be a discretization.
For $i \in [n]$ let $Z_i$ be the set of all $\vx \in M$ such that the $i$-th coordinate~$S_i$ of~$S$ takes the value~0.
We call~$S$ \emph{topologically nondegenerate} if

\smallskip

\noindent
(a) For all $i \in [n]$ both $Z_i$ and $M \backslash Z_i$ are $m$-dimensional topological manifolds.

\smallskip

\noindent
(b) For all $i \in [n]$ the boundary  $bd(Z_i)$ in~$M$ is a union of finitely many $m-1$-dimensional topological manifolds.

\smallskip

\noindent
(c) For all $i, j \in [n]$ with $i \neq j$ the intersection $bd(Z_i) \cap bd(Z_j)$ is a union of finitely many compact topological manifolds of dimensions
$\leq m-2$.

\smallskip

\noindent
(d) If $\vx \in bd(Z_i) \cap bd(Z_j) \cap int(M)$ and $U$ is a neighborhood of~$\vx$,  then for every $f: \{i, j\} \rightarrow \{0,1\}$ there exists a nonempty $V_f \subset U$ such that $S_i(\vy) = f(i)$ and $S_j(\vy) = f(j)$ for all~$\vy \in V_f$.
\end{definition}

For example, the discretizations for our flows $\bD_1(f, \vg), \bD_2(f, \vg)$ are topologically nondegenerate. In fact, whenever $M$ a product of~$N$ compact nondegenerate intervals, with
$\vx = (x_1, \ldots, x_N)$ and $S_i$ takes the value $0$ or~$1$ depending on whether~$x_i$ is below or above some threshold~$\Theta_i$ the resulting discretization~$S$ will be topologically nondegenerate.

\begin{lemma}\label{lem:nosimultaneous}
Let $\bD$ be a continuous flow and let~$M$ be a compact $m$-dimensional manifold that is forward-invariant for~$\bD$. Let $S: M \rightarrow \{0,1\}^n$ be a topologically  nondegenerate discretization of~$M$, and let the sets  $Z_i$ be as in Definition~\ref{def:BoolODedef}. Let $i \neq j$ and suppose that $\vx(0)$ is an initial condition  and
$0 < t_1 < t_3$ are times with $\{\vx(t): \ t \in [0, t_3]\}$ contained in the interior of~$M$ such that

\smallskip

\noindent
(i) $\vx(t_1) \in bd(Z_i) \cap bd(Z_j)$.

\smallskip

\noindent
(ii) For all $\vy(0)$ in some neighborhood~$U$  of $\vx(0)$ we have

$$|\{t \in [0, t_3]: \ \vy(t) \in bd(Z_i) \cap bd(Z_j)\}| \leq 1.$$

Then there exists a neighborhood $W$ of $\vx(0)$ and times $t_0, t_2$ with $0 \leq t_0 < t_1 < t_2 \leq t_3$ such that the set
\begin{equation}\label{eqn:nosimset}
NS(i, j) = \{\vy(0): \ \forall \, t \in [t_0, t_2] \ \vy(t) \notin  bd(Z_i) \cap bd(Z_j)\}
\end{equation}
contains a dense open subset~$V$ of~$W$.
\end{lemma}

\smallskip

\noindent
\textbf{Proof:}  Let  $W \subseteq U$ be a sufficiently small closed neighborhood of~$\vx(t_0)$ such that for some times $t_0 < t_1 < t_2$ trajectories that start in~$W$ don't leave~$M$ during the  time interval~$[-t_2 + t_1, 0]$.
Define a map $F: W \times [t_0, t_2] \rightarrow M \times [t_0, t_2]$ by $F(\vz(0), t) = (\vz(t - t_0), t)$. Let $K$ be the range of~$F$.

The function $F$ is continuous in both variables.  Since $W \times [t_0, t_2]$ is compact, $F$ is a homeomorphism between
$W \times [t_0, t_2]$ and $K$. Thus $K$ is a topological manifold of dimension~$m+1$.
By condition~(c) of Definition~\ref{def:BoolODedef},  the set $B = (bd(Z_i) \cap bd(Z_j)) \times [t_0, t_2]$ is a union of finitely many compact submanifolds of dimension $\leq m - 1$ of~$K$.  Consider the map $G: B \rightarrow W$ that assigns to each $(\vz, t) \in B$ the unique $\vy(0) \in W$ with $\vy(t) = \vz$.  This map is continuous by continuity of the flow, and is injective by condition~(ii).  Thus $G$ is a homeomorphism, and it follows that its range

$$R = \{\vy(0) \in W: \ \exists t \in [t_0, t_2] \ \vy(t)\in bd(Z_i) \cap bd(Z_j)\}$$

has dimension at most $m-1$. In particular, $V = W \backslash R$ is dense and open in~$W$, and the lemma follows.
 $\Box$

\begin{corollary}\label{corol:strong-implies-onestep}
Suppose $\bD$ is a flow on a compact~$m$-dimensional manifold~$M$ and $\bB = (\{0,1\}^n, f)$ is a Boolean approximation for a given topologically nondegenerate discretization $S: M \rightarrow \{0,1\}^n$.  Then strong consistency between~$\bD$ and~$\bB$ on any open~$U \subseteq int(M)$ implies that  for each $\vx \in U$ the trajectory of
$S(\vx)$ in~$\bB$ is one-stepping.
\end{corollary}

\noindent
\textbf{Proof:} Suppose we have strong consistency on $U$.  Since strong consistency implies strong consistency on each subset of~$U$, we may wlog assume that $U$ is disjoint from the boundary of $S^{-1}(\{0\})$ for all $i \in [n]$.  Assume towards a contradiction that $\vx \in U$ is such that $f_i(\vs^{\, \vx}(\tau)) \neq s_i, f_j(\vs^{\, \vx}(\tau)) \neq s_j$ for some $\tau$ and $i \neq j$.  Assume wlog $\tau = 0$.  Then we find $t_1 > 0$ so that~(i) of Lemma~\ref{lem:nosimultaneous} holds.  Moreover, by continuity and the transversality property, we can choose $t_3$ in such a way that point~(ii) of the lemma holds in some neighborhood of~$\vx$; let us for simplicity assume that~$U$ itself has this property.   But then for all $\vy \in NS(i, j) \cap U$ we must have either $f_i(\vs^{\, \vy} (0)) = s^{\vy}_i$ or $f_j(\vs^{\, \vy} (0)) = s^{\vy}_j$; otherwise
 strong consistency along the trajectory of~$\vy$ would be violated. Thus for all $\vy \in NS(i, j) \cap U$ we must have  $S(\vy) \neq S(\vx)$.  But since
 $NS(i, j)$ contains a dense open subset of some neighborhood~$W$ of~$\vx$ by Lemma~\ref{lem:nosimultaneous}, we get a contradiction with point~(d) of Definition~\ref{def:BoolODedef}.
 $\Box$

\bigskip

\begin{definition}\label{def:mono-stepping}
Let $\bB = (\{0,1\}^n, f)$ and let $\vs, \vsst \in \{0,1\}^n$.

\smallskip

\noindent
(i) We say that $\vsst$ is \emph{strictly between $\vs$ and $f(\vs)$} and write $\vs \preceq \vsst \prec f(\vs)$ if
$s_i = s_i^*$ for all $i$ with $s_i = f(\vs)_i$ and $\vsst \neq f(\vs)$.

\smallskip

\noindent
(ii) We say that \emph{the trajectory $\{\vs(\tau)\}_{\tau = 0}^\infty$ of $\vs =\vs(0)$ is monotone-stepping} if the following holds:

\begin{equation}\label{eqn:mono-stepping-trajs}
\forall \tau \ \forall \vsst  \ (\vs(\tau) \preceq \vsst \prec \vs(\tau+1)  \Rightarrow f(\vsst) = \vs(\tau +1)).
\end{equation}

\noindent
(iii) We say that \emph{$\bB$ is monotone-stepping} if all its trajectories are monotone-stepping.
\end{definition}

Clearly, if $H(\vs, f(\vs)) = 1$, then the only $\vsst$ with $\vs \preceq \vsst \prec f(\vs)$ is $\vs$ itself.  Thus all one-stepping Boolean trajectories and Boolean systems are monotone-stepping, but not \emph{vice versa.}

\begin{example}\label{ex:mono-stepping}
The function $f: \{0,1\}^3 \rightarrow \{0,1\}^3$ given by

\begin{equation}\label{eqn:mono-stepping-example}
\begin{split}
f(000) &= 110 = f(010) = f(100);  f(110) = 111;\\
f(111) &= 101; f(101) = 001; f(001) = 011; f(011) = 111
\end{split}
\end{equation}
defines a monotone-stepping Boolean system that is not one-stepping.
\end{example}

\section{Consistency and Strong Consistency: Some Basic Examples}\label{sec:examples}

In this section we include some examples that will illustrate the constructions and definitions given in previous sections and will put our main results into context.  More details and additional examples can be found in~\cite{LowDimEx16}.

\subsection{$\bD_1(f, \gamma)$ for Boolean constants}\label{subsec:constsec1}

Let $f: \{0,1\} \rightarrow \{0,1\}$ be a Boolean constant, that is, $f(s) \equiv 0$ or $f(s) \equiv 1$. Since we are working in one dimension here, the ODE that defines~$\bD_1(f, \vg)$ becomes

\begin{equation}\label{eqn:standardform-1dim}
\dot{x}=\gamma (g(x)+ 6 Q(\vx)).
\end{equation}

If $Q$ takes only values in the interval $[0, 1/6)$ (for  $f(s) \equiv 0$) or $(5/6, 1]$ (for $f(s) \equiv 1$), then~(\ref{eqn:standardform-1dim}) has only one globally asymptotically stable equilibrium that all trajectories will approach monotonically, and strong consistency on the whole state space follows. This will clearly be the case for the conversion scheme~$\cQ_W$, which returns a constant~$Q$.

Other conversion schemes may not be so well-behaved. For example, the Boolean expression $(s \wedge \neg s)$ also represents the constant function  $f(s) \equiv 0$,  but gets translated by $\cQ_{Rc}, \cQ_{Rd}$ into $Q$ that takes all values on the interval~$[0, 1/4]$. These conversion methods return the following instance
of~(\ref{eqn:standardform-1dim}):

\begin{equation}\label{eqn:swedgenotseq}
\dot{x}= \gamma (g(x)+ 6  L(x)(1-L(x))).
\end{equation}

\begin{figure}[h!]\caption{$\dot{x}= g(x)+ 6  L(x)(1-L(x))$}
	\centering
	\includegraphics[width=3 in, height=3 in]{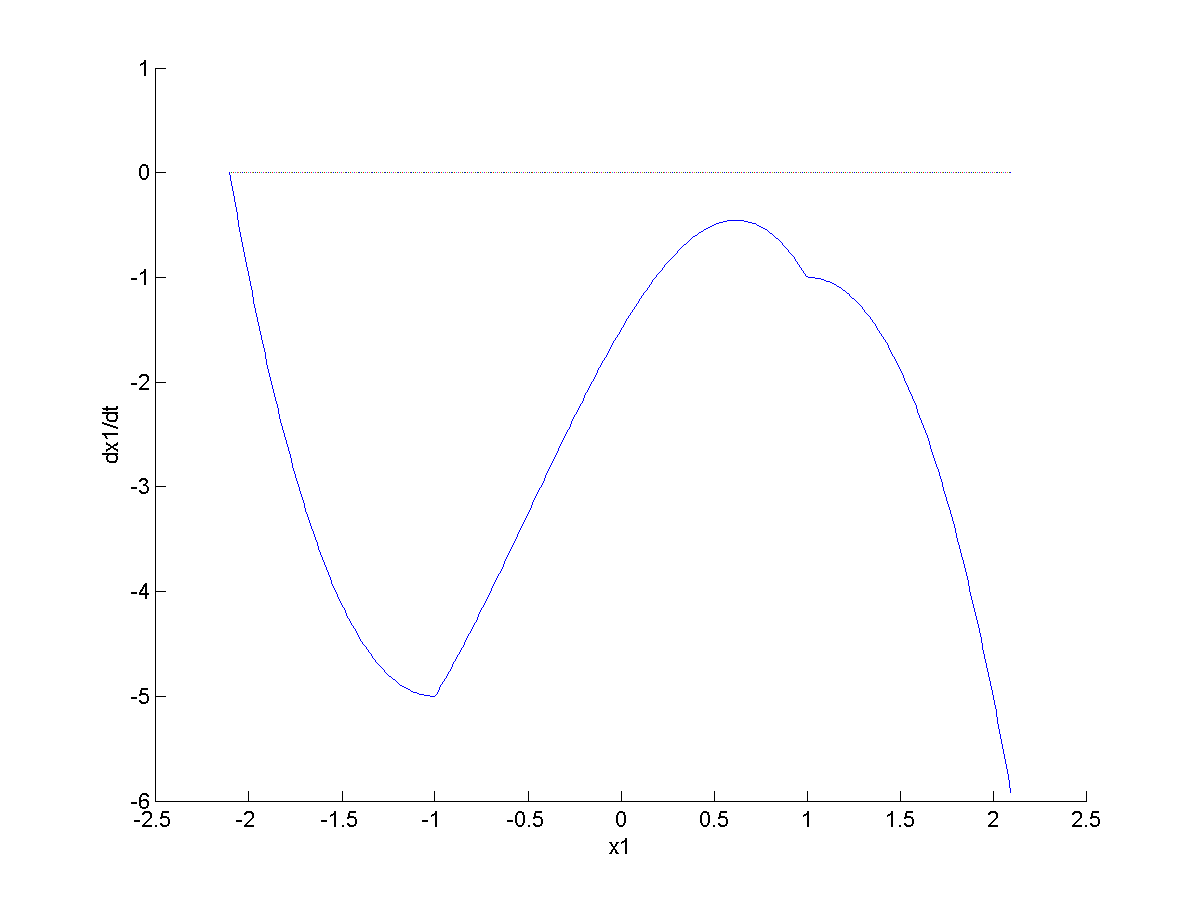}
	\label{fig:cubic}
\end{figure}

Figure~\ref{fig:cubic} shows the right-hand side for $\gamma = 1$, and it can be seen that strong consistency will still hold.

However, more complicated Boolean expressions for a contradiction may cause problems. For example, consider~$f$ that is given by the Boolean expression~$f(r) = \sigma \wedge \neg \sigma$, where $\sigma =s \wedge s \wedge s \wedge s$.  This may be recursively converted into the following instance of~(\ref{eqn:standardform-1dim}):

\begin{equation}\label{eqn:k4eqn}
\dot{x}= \gamma(g(x)+6L(x)^4(1-L(x)^4)).
\end{equation}

\begin{figure}
\caption{$\dot{x}=g(x)+6L(x)^4(1-L(x)^4)$}
	\centering
		\includegraphics[width=3 in, height=3 in]{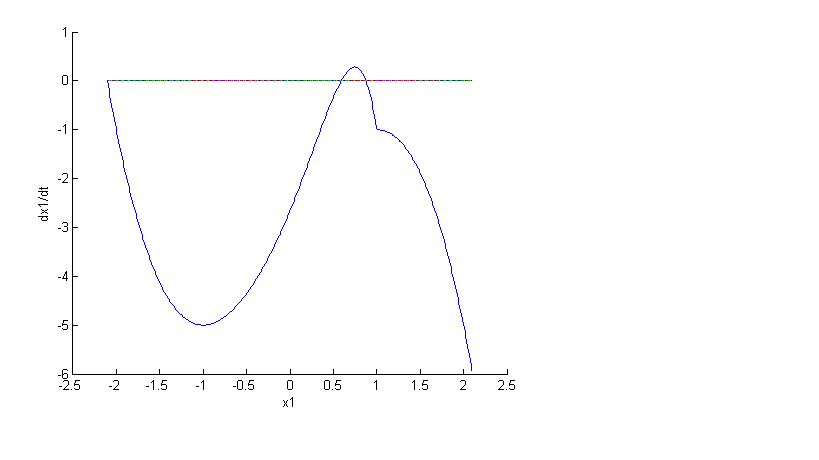}
	\label{fig:problemcubic}
\end{figure}

Figure~\ref{fig:problemcubic} shows the graph for~$\gamma = 1$.

The system has three fixed points $r_1=x^-, r_2=.58875, r_3=.87703$. We get strong consistency on $[x^-, r_2)$ (which is a universal set), but not on the whole state space.

In view of these potential problems one may wonder why we consider conversion methods other than~$\cQ_W$ at all.  There are at least three reasons for this.
As we already mentioned in Subsection~\ref{subsec:convert-schemes}, speed of computing the conversion is one of them.  Second, it seems a good idea to allow for more generality at little extra cost of the exposition, especially since some of our recursive conversion schemes are quite natural for simple updating functions.
Third, suppose for example that $f_i = s_1 \wedge s_3$.  Even~$\cQ_W$ will translate this into a composition of a quadratic polynomial with~$L$.  However, if we investigate the behavior along a trajectory for which which $s_3 \equiv 0$, then $f_i$ will behave along this trajectory as a Boolean constant in exactly the same way as any contradiction, while $Q_i$ is not a constant. We want to build up a set of general tools that can be applied in such situations.

\subsection{Conversions of copy-negation}\label{subsec:copy-negation}

There are two nonconstant Boolean functions~$f:\{0,1\} \rightarrow \{0,1\}$ of one variable:  $f(s) = s$ (the ``copy'' function), and
$f(s) = \neg s$ (``copy-negation''). For the latter function all Boolean trajectories satisfy

$$\dots \mapsto 0 \mapsto 1 \mapsto 0 \mapsto 1 \mapsto \dots$$

In order to have consistency  with this function we need some type of oscillations, which is not possible in a one-dimensional flow, at least if the discretization is based on a single threshold. In particular, we cannot have consistency along any trajectory of $\bD_1(f, \vg, \cQ)$, for any choice of $\vg$ and~$\cQ$.

Notice that the problem here is caused because $f$ copies the negation of its variable to the variable \emph{itself,} which is similar to the problem caused by \emph{self-regulation} in Glass networks (see Subsection~\ref{subsec:related-work}).  Let us try to sidestep the problem by considering a two-dimensional system where the copying and negating steps are separated.

Let $f = (f_1, f_2): \{0,1\}^2 \rightarrow \{0,1\}^2$ be defined by:

\begin{equation}\label{eqn:2-regs}
f_1(s) = \neg s_2 \qquad  f_2(s) = s_1.
\end{equation}

This system has a single orbit  given by

\begin{equation}\label{eqn:fcncycle}
00 \mapsto 10 \mapsto 11 \mapsto 01 \mapsto 00.
\end{equation}

Define $\bD = \bD_1(f, (1,1))$ by any of the conversion schemes $\cQ_W, \cQ_{Rc}, \cQ_{Rd}, \cQ_{RF}$.  This gives:

\begin{equation}\label{eqn:copynegDE}
\begin{split}
\dot{x}_1 &= g(x_1)+ 6  (1-L(x_2)),\\
\dot{x}_2 &= g(x_2)+ 6  L(x_1).
\end{split}
\end{equation}

\begin{figure}[ht]
\centerline{\includegraphics[width=5in]{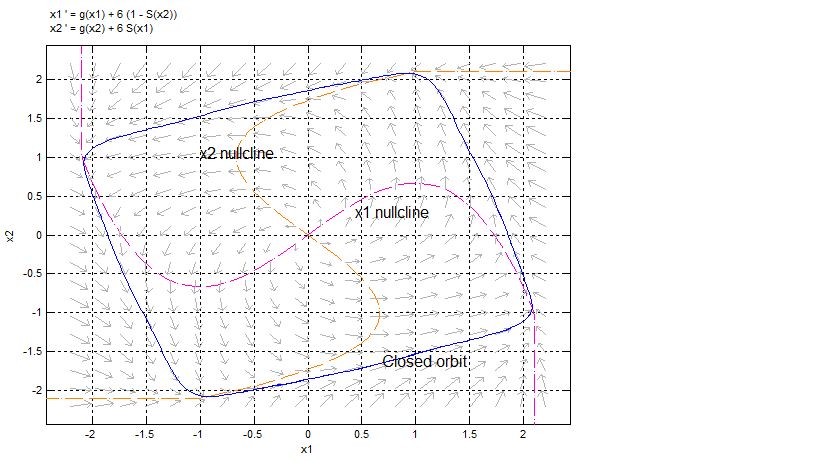} }
\caption{Nullclines and direction arrows for $\bD_1(f, (1, 1))$.}
\label{F1}
\end{figure}

Since $n = 2$, we have the luxury of being able to perform an easy phase-plane analysis of $\bD$.
Figure~\ref{F1} gives the phase portrait.

The two nullclines intersect at~$(0,0)$, which is the only steady state. The eigenvalues of the Jacobian at the steady state have positive real parts, and the system has a unique limit cycle that will be approached by all trajectories that start off the origin.  It can easily be seen by inspection of~Figure~\ref{F1}  that the flow~$\bD$ is strongly consistent with~$f$ on the set~$U = [x^-, x^+]^2 \backslash \{0, 0\}$. This will remain true even if we allow arbitrary vectors $\vg$; only the shape of the limit cycle and nullclines in Figure~\ref{F1} will change for $\gamma_1 \neq \gamma_2$, but not the qualitative behavior of the system.

Recall that for one-dimensional flows we cannot get any consistency with copy-negation whatsoever.   In a sense, we added just one dummy variable, and we got  as much consistency with the Boolean dynamics as one could possibly hope for. Note that  $\bD$ is really nothing else but $\bD_2(f, (1,1))$ with the roles of variables reversed, and the consistency result obtained here carries over to $\bD_2(f, (1,1))$ for one-dimensional copy-negation~$f$.   This may exemplify a more general tendency of  intermediate variables to favor consistency.  Our main result in Section~\ref{sec:Consistency-Stepping} also points in the same direction.

\subsection{An example that is not monotone-stepping}\label{subsec:quasipersec}

Define a Boolean updating function $f: \{0,1\}^4 \rightarrow \{0,1\}^4$  by

\begin{equation}\label{eqn:22-regs}
\begin{split}
f_1(s) = \neg s_2 &\qquad  f_2(s) = s_1,\\
f_3(s) = \neg s_4 &\qquad  f_4(s) = s_3.
\end{split}
\end{equation}

This system is the direct product of the Boolean
system defined by~(\ref{eqn:2-regs}) of Subsection~(\ref{subsec:copy-negation}) with itself.  There are four pairwise disjoint orbits of length four each in this system:

\begin{equation}\label{eqn:22-trajs}
\begin{split}
0000 \mapsto 1010 \mapsto 1111 \mapsto 0101 \mapsto 0000,\\
0010 \mapsto 1011 \mapsto 1101 \mapsto 0100 \mapsto 0010,\\
0011 \mapsto 1001 \mapsto 1100 \mapsto 0110 \mapsto 0011,\\
0001 \mapsto 1000 \mapsto 1110 \mapsto 0111 \mapsto 0001,
\end{split}
\end{equation}
and it can be seen that the system is not monotone-stepping. Under several of our conversion schemes the
flow~$\bD_1(f, \vec{1})$ will be generated by the DEs

\begin{equation}\label{eqn:copynegDE4}
\begin{split}
\dot{x}_1 &=\gamma_1 ( g(x_1)+ 6  (1-L(x_2))),\\
\dot{x}_2 &= \gamma_2 (g(x_2)+ 6  L(x_1)),\\
\dot{x}_3 &= \gamma_3 (g(x_3)+ 6  (1-L(x_4))),\\
\dot{x}_4 &= \gamma_4 (g(x_4)+ 6  L(x_3)),
\end{split}
\end{equation}
with no interaction between the variables $(x_1, x_2)$ and $(x_3, x_4)$.  Thus~(\ref{eqn:copynegDE4}) can be treated as a direct product of two flows to which the results of Subsection~\ref{subsec:copy-negation} apply.

It is not hard to infer that all trajectories of~(\ref{eqn:copynegDE4}) that start in the set

$$U = [x^-, x^+]^4 \backslash \{\vx: \ x_1 = x_2 = 0 \ \vee \ x_3 = x_4\}$$
are consistent with~$f$, regardless of the choice of~$\gamma$.  Strong consistency on open sets in this example is precluded by Corollary~\ref{corol:strong-implies-onestep}.  For generic choices of $\vg$ the ODE dynamics on the attractor will become aperiodic (but not chaotic, only quasiperiodic) which  precludes strong consistency with~$f$ for even a single trajectory of~(\ref{eqn:copynegDE4}).

After changing the order of variables in~(\ref{eqn:copynegDE4}) the system becomes $\bD_2(f, \vg)$ for
$f: \{0,1\}^2 \rightarrow \{0,1\}^2$ defined by

\begin{equation}\label{eqn:cn-22-regs}
f_1(s) = \neg s_1 \qquad  f_2(s) = \neg s_2,
\end{equation}
and analogous observations about consistency and strong consistency apply.  The choice of the word ``analogous'' rather than ``the same'' in the preceding sentence was very deliberate, since when treating the same flow as $\bD_2(\cdot,\vg)$ rather than $\bD_1(\cdot,\vg)$, we are comparing its  dynamics with a different Boolean system!

Note that the function~$f$ defined by~(\ref{eqn:cn-22-regs}) has two disjoint periodic orbits

\begin{equation}\label{eqn:cn-22-trajs}
00 \mapsto 11 \mapsto  00, \qquad
01 \mapsto 10 \mapsto 01,
\end{equation}
and, similar to the function~(\ref{eqn:22-regs}), is not monotone-stepping.

\section{Consistency for Monotone-Stepping Boolean Systems}\label{sec:Consistency-Stepping}

In this section we will prove the following result.

\begin{theorem}\label{thm:mainth}
Let $f: \{0,1\}^n \rightarrow \{0,1\}^n$ be a Boolean function, let~$\cQ$ be a conversion scheme, and let
$\vg^{\, -} = (\gamma_1, \ldots , \gamma_n)$ be a fixed vector of positive reals.  Then there exist $\mu > 0$ and nonempty open~$U^s \subset W^s$ for $s \in \{0,1\}^n$ such that for every extension of~$\vg^{\, -}$ to a $2n$-dimensional vector $\vg$ of positive reals with $\gamma_{i+n} < \mu$ for all~$i \in [n]$ the following holds in~$\bD_2(f, \vg, \cQ)$:

\smallskip

\noindent
(i) For every $\vs \in\{0,1\}^n$ whose trajectory in $\bB = (\{0,1\}^n, f)$ is monotone-stepping the flow $\bD_2(f, \vg, \cQ)$ is consistent with~$\bB$ on~$U^s$.

\smallskip

\noindent
(ii) For every $\vs \in\{0,1\}^n$ whose trajectory in $\bB = (\{0,1\}^n, f)$ is one-stepping the flow  $\bD_2(f, \vg, \cQ)$ is strongly consistent with~$\bB$.

\smallskip

\noindent
(iii) If $\bB = (\{0,1\}^n, f)$ is monotone-stepping, then the flow $\bD_2(f, \vg, \cQ)$ is consistent with~$\bB$.

\smallskip

\noindent
(iv) If $\bB = (\{0,1\}^n, f)$ is one-stepping, then the flow $\bD_2(f, \vg, \cQ)$ is strongly consistent with~$\bB$.
\end{theorem}

Before proving the theorem, let us make some comments.  By choosing $\mu$ sufficiently small relative to $\min \vg^{\, -}$ we can assure that the values of the signaling variables $x_{i+n}$ change at a slower timescale than the values of the signature variables~$x_i$.  Thus Theorem~\ref{thm:mainth} tells us that the relevant consistency results hold, with fixed witnesses $U^s$ or $U = \bigcup \{U^s: \ s \in\{0,1\}^n\}$, \emph{under sufficient separation of time scales.}

By our definitions, point~(iii) of the theorem follows from point~(i), and point~(iv) follows from point~(ii).  Moreover, point~(ii) follows from point~(i) since by Proposition~\ref{prop:mono+cons-implies-strong} for one-stepping trajectories strong consistency is the same as consistency.

\bigskip

\noindent
\textbf{Proof of Theorem~\ref{thm:mainth}:} It suffices to prove part~(i).

First let us introduce an additional parameter~$\delta$ in  the definition of the sets~$W^s_i$ of~(\ref{eqn:Wsi-def}):

\begin{equation}\label{eqn:Wsidelta-def}
\begin{split}
W^s_i(\delta) &= \{\vx: \ x_i < -1 + \delta\} \ \mbox{if} \ s_i = 0,\\
W^s_i(\delta) &= \{\vx: \ x_i > 1- \delta\} \ \mbox{if} \ s_i = 1.
\end{split}
\end{equation}

Let $W^s(\delta) = \bigcap_{i \in [n]} W^s_i(\delta)$.  Since the functions $Q_j$ are continuous, by~(\ref{eqn:Q-reflects-Boolean}) there exists $\delta > 0$ such that for all~$s$ and $i \in [n]$ and all fixed $x_1, \dots , x_n \in W^s(\delta)$ the DE for ${x}_{i+n}$ in~(\ref{eqn:standardform2}) has a globally attracting equilibrium.  Let us fix a sufficiently small such~$\delta$ throughout this proof.

Now there exist positive constants~$\beta, \alpha > 0$ that depend only on~$\delta$ such that for all $i \in [n]$ and states~$\vx$ that satisfy either

\begin{equation}\label{eqn:xi-past-delta1}
x_{i} \in [-1 + \delta, 1] \quad \mbox{and} \quad x_{i+n} \geq 2/3 - \alpha
\end{equation}

or

\begin{equation}\label{eqn:xi-past-delta2}
x_{i} \in [-1, 1 - \delta] \quad \mbox{and} \quad x_{i+n} \leq -2/3 + \alpha
\end{equation}

we have

\begin{equation}\label{eqn:minspeedxi}
|\dot{x}_{i}| > \gamma_i\beta.
\end{equation}

To see this, first note that $x_{i+n} \geq 2/3 - \alpha$ iff $L(x_{i+n}) \geq 5/6 - \alpha/2$ and $x_{i+n} \leq -2/3 + \alpha$ iff $L(x_{i+n}) \leq 1/6 + \alpha/2$. Choose $\alpha$ with $0 < \alpha < \delta$ such that the unstable equilibrium $x^\circ$ of
$\dot{x} = g(x) + 5 + 3\alpha$ satisfies $x^\circ < -1 + \delta$. Notice that  for this choice of $\alpha$  the unstable equilibrium $x^\circ$ of
$\dot{x} = g(x) + 5 - 3\alpha$ satisfies $x^\circ > 1 - \delta$, and hence $\dot{x}_{i} \neq 0$ on the compact set of all those $\vx$ that satisfy~(\ref{eqn:xi-past-delta1}) or~(\ref{eqn:xi-past-delta2}). Now the existence of~$\beta$ follows from the Extreme Value Theorem.  We will fix~$\alpha, \beta$ as above for the remainder of this proof.  For technical reasons we will assume $\alpha < 2/15$.

Now let us define analogues of the sets~$W^s_i(\delta)$ for the variables $x_{i+n}$.  For $\eps \geq 0$ let

\begin{equation}\label{eqn:Wsidelta-def}
\begin{split}
V^s_i(\eps) &= \{\vx: \ x_i < -2/3 + \eps\} \ \mbox{if} \ s_i = 0,\\
V^s_i(\eps) &= \{\vx: \ x_i > 2/3 - \eps\} \ \mbox{if} \ s_i = 1.
\end{split}
\end{equation}

In analogy with~$W^s_i$, we will write $V^s_i$ instead of~$V^s_i(0)$.

For each $s \in \{0,1\}^n$ define the following sets

$$\Delta(s) = \{ i \in [n]: \, f_i(s) \neq s_i\}, \qquad \Gamma(s) = [n]\backslash \Delta(s).$$

Let

$$\alpha_s = \frac{n - |\Delta(s)|}{n}\alpha.$$

Define $R^s \subset M$ as the open set of all states in the set

$$PR^s = \bigcap_{i \in \Delta(s)} W^s_i(\delta) \times \bigcap_{i \in \Gamma(s)} W^s_i \times \bigcap_{i \in \Gamma(s)} V^s_i(\alpha_s)$$

that satisfy the following condition for all $i \in [n]$:

\begin{equation}\label{eqn:pre-release}
\vx \in W^s_i(\delta) \backslash W^s_i \quad \mbox{implies} \quad \vx \in V^{f(s)}_i(\alpha_s).
\end{equation}

Let

$$U^s = R^s \cap W^s \cap \prod_{i \in \Gamma(s)} V^s_i.$$

Then $U^s$ is nonempty, open, and $S(\vx) = s$ for every~$\vx \in U^s$.

For a given state $\vx \in M$, a given~$s \in \{0,1\}^n$, and $j \in [2n]$ consider the DE for $x_j$ in~(\ref{eqn:standardform2}).  Let $x^*_j(\vx, s)$ denote the locally stable
equilibrium that is $< -1$ if $s_j = 0$ or $s_{j -n} = 0$ and this equilibrium exists, or the locally stable
equilibrium that is $> 1$ if $s_j = 1$ or $s_{j -n} = 1$ and this equilibrium exists.  Similarly, let $x^\circ_j(\vx, s)$ denote the unstable equilibrium in $(-1, 1)$ if it exists. In our arguments, we will omit the parameters $\vx, s$ if they are implied by the context.

Equipped with this terminology,  let us consider as a warm-up what happens to an ODE trajectory that starts in $U^s$ when $s$ is a fixed point of~$f$. Then $\Delta(s) = \emptyset$ and $\Gamma(s) = [n]$. In this case the equilibria $x^\circ_i$ do not exist, and $x^*_i$ will be globally attracting in the dynamics of~$x_i$ for all~$i \in [n]$. Thus the trajectory will never leave $W^s$, which implies that the variables $x_{i+n}$ will keep moving away monotonically from the interval $[-2/3, 2/3]$.  It follows that~$U^s$ is forward invariant in this case; in particular, no change in the Boolean state will occur at any future time, and we get strong consistency on~$U^s$ with the Boolean trajectory of~$s$.

If $s$ is not a fixed point of~$f$, things become more complicated. Let us call a Boolean state~$s$ \emph{monotone for~$f$} if $f(s) = f(s')$ for all $s'$ with $s \preceq s' \prec f(s)$. Clearly, a monotone-stepping trajectory is one that visits only monotone states.  In view of what we have already shown for fixed points, if suffices to prove that
there exist $\mu > 0$ and $T_{min} > 0$  such that for every extension of~$\vg^{\, -}$ to a $2n$-dimensional vector $\vg$ of positive reals with $\gamma_{i+n} < \mu$ for all~$i \in [n]$ and every
$s$ that is monotone but not a fixed point of~$f$ and every $\vx(0) \in U^s$  the following conditions hold:

\begin{itemize}
\item[(A)] $\vx(T) \in U^{f(s)}$,
\item[(B)] $x_i(t) \neq 0$ for every $i \in \Gamma(s)$ and all $t \in [0, T]$,
\item[(C)] For every $i \in \Delta(s)$ there exists exactly one $t \in [0, T]$ with $x_i(t) = 0$.
\end{itemize}

Let us first make a couple of remarks.

The requirement that the times~$T$ are bounded away from~$0$ by a fixed $T_{min}$ is needed to ensure condition~(e) of the transversality property.  But the existence of such~$T_{min}$ will follow automatically from~(A), because as $s$ is assumed not to be a fixed point, at least one of the $x_i$'s has to traverse the entire interval~$[-1, 1]$ between times~$0$ and~$T$, and since  $|\dot{x}_i| \leq 5\gamma_i$ (see~(\ref{eqn:derivest-upper}) below), this will take at least $\frac{2}{5\gamma_i}$ time units. Taking $\gamma_i = \max \vg^{\, -}$ here gives a required $T_{min}$.

Since there are only finitely many Boolean states to consider, we can simplify our task by determining
a suitable $\mu = \mu(s)$ for each relevant~$s$ separately and then taking~$\mu$ as the minimum.
Thus for the remainder of the proof we will fix $s$ that is monotone and assume that $\vx(0) \in U^s$.

We will argue that under sufficient separation of time scales there exists $T$ such that (A)--(C) are satisfied and the trajectory of~$\vx(0)$ on the time interval
$[0, T]$ can be partitioned into three types of \emph{episodes.}

For all times~$t$ in a \emph{Type~I episode} we will have $\vx \in R^{s'}$ for some $s'$ with
$s \preceq s' \prec f(s)$.  Notice that such~$s'$ must be monotone and must satisfy $f(s') = f(s)$, $\Delta(s') \subseteq \Delta(s)$ and  $\Gamma(s') \supseteq \Gamma(s)$.  Clearly, time~$t = 0$ belongs to a Type~I episode with $s' = s$.

During a Type~I episode, all $x_i$ for $i \in \Gamma(s')$ move towards the equilibrium $x^*_i$ which is globally stable in this case.
Similarly, by the choice of~$\delta$, each variable  $x_{i+n}$ for $i \in \Gamma(s')$ will move towards the unique globally stable equilibrium~$x^*_{i+n}$.
In particular, throughout a
Type~1 episode, the system will stay within the set
$$\bigcap_{i \in \Gamma(s')} W^{s'}_i \cap \bigcap_{i \in \Gamma(s')} V^{s'}_i(\alpha_{s'}) \subseteq \bigcap_{i \in \Gamma(s)} W^{s}_i \cap \bigcap_{i \in \Gamma(s)} V^{s}_i(\alpha_{s'}).$$

In contrast, each variable  $x_{i+n}$ for $i \in \Gamma(s')$ will move away towards the unique globally stable equilibrium in $V^{f(s')}_i = V^{f(s)}_i$.
This movement will continue until either the variable enters $V^{f(s)}_i$, or some $x_i$ for $i \in \Delta(s') \subseteq \Delta(s)$ enters the interval $[-1 + \delta, 1 - \delta]$, that is, $\vx$ leaves $W^{s'}_i(\delta)$.  In the former case, $x^*(f(s))$ becomes the globally stable equilibrium of the DE~(\ref{eqn:standardform2}) for $x_i$, and $x_i$ starts moving towards the boundary of $[-1+\delta, 1 - \delta]$.  This movement can be arrested only if $x_{i+n}$ changes direction, which in turn requires some variable $x_{i'}$ to enter the interval $[-1 + \delta, 1 - \delta]$.  The upshot is that \emph{a Type~I episode ends at a time $t_i$ when $x_i(t_{s'}) \in \{-1 + \delta, 1 - \delta\}$ for at least one variable $x_i$ with $i \in \Delta(s') \subseteq \Delta(s)$.}

The last condition marks the onset of a \emph{Type~II episode for~$i$.}  We can see from the above discussion that Type~II episodes for different~$i \in \Delta(s)$ may overlap or even occur simultaneously, but we will see that only one such episode for each~$i \in \Delta(s)$ can occur in the interval~$[0, T]$.

We want to show that under sufficient separation of time scales, during a Type~II episode for~$i$, the variable~$x_i$ will move monotonically from~$-1+\delta$ into the interval~$(1, x^+]$ or from~$1-\delta$ into the interval~$[x^-, -1)$ while not much movement of variables $x_{j+n}$ occurs. For ease of exposition, consider the onset of a Type~II episode for~$i$ at time~$t_i$ with $x_i(t_i) = -1 + \delta$; the case $x_i(t_i) = 1 - \delta$ is symmetric.  By~(\ref{eqn:minspeedxi}), the variable $x_{i}$ will keep monotonically increasing and will reach the interval $(1, x^+]$ at some time $T_i > t_i$ with

\begin{equation}\label{eqn:xi-reachtime}
T_i - t_i < \frac{1}{\gamma_{i}\beta}
\end{equation}
unless $x_{i+n}(t) < 2/3 - \alpha$ for some time $t > t_i$ with $t - t_i < \frac{1}{\gamma_{i^*}\beta}$.

Note that since $Q_i$ is assumed to take values only in the interval~$[0,1]$, the shape of~$g$ implies that we have at all times~$t$ and for all $j \in [2n]$

\begin{equation}\label{eqn:derivest-upper}
 |\dot{x}_{j}(t)| \leq 5\gamma_{j}.
\end{equation}

Now consider~$t_{i_0}$ that marks the onset of the first of~$k$ overlapping Type~II episodes for $i \in I \subseteq \Delta(s')$.
Then we must have $k \leq |\Delta(s')| \leq |\Delta(s)|$, and the combined length of these Type~II episodes is at most~$k\min \vg^{\, -}$.  Let

 \begin{equation}\label{eqn:mufirst}
0 < \mu < \frac{\beta\min \vg^{\, -}\alpha}{5n},
\end{equation}
and assume $\gamma_{j+n} < \mu$ for all $j \in [n]$.

Condition~(\ref{eqn:pre-release}) applies to all $i \in I$ at time $t_0$, and for each such~$i$ we must have $\vx(t) \in V^{s'}_i(\alpha_{s'})$ at some time~$t$ in these overlapping Type~II episodes. But by~(\ref{eqn:derivest-upper}), for each $i \in I$ it will take $x_{i+n}$ more than $\frac{k}{\beta\min \vg^{\, -}}$ time units to travel a distance of more than $\frac{k\alpha}{n}$ units.  This shows that all $x_i$ for $i \in I$ will reach their destinations in
$[x^-, -1) \cup (1, x^+]$ before the system has a chance to move out of $\bigcap_{i \in I} V^{s'}_i(\alpha)$.  Similar considerations apply to the variables $x_{j+n}$ for $j \in \Gamma(s')$: while we don't have not much control over~$Q$ during a Type~II episode and these variables may move briefly in the wrong direction, they cannot get very far; each of them will move by less than $\frac{k\alpha}{n}$ units.

Let us make the observation more precise in the following way: Consider Type~II episodes as above and let $s''$ be such that $s''_i = s'_i$ for all~$i \in [n] \backslash I$ and $s''_i = 1 - s'_i$ for all~$i \in I$.  Then $s' \prec s'' \preceq f(s)$.
If $s'' \neq f(s)$, then $s''$ is monotone and $\Delta(s'') = \Delta(s') \backslash I$, $\Gamma(s'') = \Gamma(s') \cup I$. In this case,
after completion of these Type~II episodes, the system will have entered a state in~$R^{s''}$, and the next Type~I episode starts.

This alternation of Type~I and Type~II episodes will continue until all variables $x_i$ with $i \in \Delta(s)$ have crossed the interval $[-1, 1]$, exactly once, assuming the value~$0$ at exactly one time along the way, which implies~(C).  Moreover, all variables $x_i$ with~$i \in \Gamma(s)$ will stay in~$[x^-, -1) \cup  (1, x^+]$, which implies~(B). Eventually, after completing the last Type~II episode, the system will reach a state $\vx(T^-) \in R^{f(s)} \cap W^{s''}$ where $s'' = f(s)$.

Such time $T_i$ will mark the onset of a \emph{Type~III episode.}  We need only to show that under sufficient separation of time scales during a Type~III episode the trajectory will move from~$R^{f(s)}$ into $U^{f(s)}$ while staying inside $W^{f(s)}$.  This will imply~(A) without violating~(B) or~(C) and complete the proof of the theorem.

We will not need to assume that $f(s)$ is monotone, but we still need to consider the set $\Gamma(f(s))$ and show that at some time $T \geq T^-$ we have
$\vx(T) \in \bigcap_{i \in \Gamma(f(s))} V^{f(s)}_i$ while $\vx(t) \in W^{f(s)}$ for all $t \in [T^-, T]$. Notice that at time $T^-$ for each $i \in [n]$ the variables
$x_i, x_{i+n}$ will be on the same side of~$0$ and head in the direction of the locally stable equilibria $x^*_i$ and $x^*_{i+n}$. If $x_{i+n}(T^-) \in [-2/3, 2/3]$ then $x_{i+n}$ must have entered this interval during a Type~II episode (after leaving it and releasing~$x_i$ if $i \in \Delta(s)$), and can have subsequently moved deeper into this interval only during some of the Type~II episodes when we didn't have control over~$Q$.  Thus by~(\ref{eqn:derivest-upper}) the distance of $x_{i+n}$ from the nearest endpoint of this interval can be at most $5\gamma_{i+n}K$, where $K$ denotes the combined  duration of all Type~II episodes since time $t=0$.  By~(\ref{eqn:xi-reachtime}), we can estimate

\begin{equation}\label{eqn:Lestim}
K \leq \frac{n}{\beta \min \vg^{\, -}}.
\end{equation}

On the other hand, while $\vx(t) \in W^{f(s)}$,
by~(\ref{eqn:Q-reflects-Boolean}) and the shape of~$g$ we have

\begin{equation}\label{eqn:derivest-lower}
\gamma_{i+n} \leq |\dot{x}_{i+n}(t)|.
\end{equation}

Let $T = T^- + \frac{5n}{\beta \min \vg^{\, -}}$.
Now~(\ref{eqn:Lestim}) and~(\ref{eqn:derivest-lower}) imply that the trajectory of $\vx(T^-)$ will reach $V^{f(s)}_i$ at or before time~$T$ as long as it stays inside $W^{f(s)}$ during the time interval~$[T^-, T]$.  But the trajectory can leave~$W^{f(s)}$ only after some $x_{j+n}$ has crossed over to the other side of the interval
$[-2/3 + \alpha, 2/3 - \alpha]$, which in view of~(\ref{eqn:derivest-upper}) takes at least~$\frac{4/3 - 2\alpha}{5 \gamma_{j+n}}$ time units. Since we assumed that $\alpha < 2/15$, we have $\alpha < \frac{4/3 - 2\alpha}{5}$, and it follows from the choice of~$\mu$ in~(\ref{eqn:mufirst}) that $\vx(t)$ will stay inside $W^{f(s)}$ during the time interval~$[T^-, T]$ as long as $\gamma_{i+n} < \mu$ for all~$i \in [n]$. $\Box$

\bigskip

Let us remark that we may see consistency even in systems that are not monotone-stepping;  Subsection~\ref{subsec:quasipersec} gives an example.  In this example, the flow $\bD_2(f,\vg)$ is a direct product of two flows with one-stepping Boolean approximations whose variables don't interact. Theorem~\ref{thm:mainth} clearly applies to to each factor and thus generalizes to direct products of monotone-stepping systems. It remains an open problem to find the most general assumptions under which the conclusion of Theorem~\ref{thm:mainth} holds.

\section{Discussion}\label{sec:discussion}

\subsection{Related results}\label{subsec:related-work}

The literature contains a few results on consistency between differentiable flows and Boolean systems; here we discuss how these are related to our Theorem~\ref{thm:mainth}.  The best-known examples are the results for Glass networks. These models were developed by Glass in~\cite{Glass1}--\cite{Glass4} based on earlier work on similar models
by Glass and Kauffman~\cite{GlassKauffman}.  The paper~\cite{Edwards} gives a review of the major results and a unifying approach to the extensive literature on these networks.
For suitable choices of the functions~$F_i$ and decay constants $\lambda_i$ in~(\ref{eqn:GlassODEs}), the resulting flows are consistent with their Boolean approximations. Consistency becomes problematic though
in the presence of \emph{self-regulation,} that is, if $F_i$ depends on~$x_i$.
As our work in Subsection~\ref{subsec:copy-negation} shows, similar problems may occur in our systems~$\bD_1(f, \vg)$. While the results about consistency in Glass networks do not require an assumption of monotone-stepping trajectories, Corollary~\ref{corol:strong-implies-onestep} applies and strong consistency does require one-stepping Boolean trajectories.

In contrast, Terman \emph{et al.}  \cite{TAWJ} constructed a class of excitatory-inhibitory neuronal networks whose
ODE models have natural Boolean approximations and showed that \emph{every} Boolean system can be translated into a model in this class so that we have a form of strong consistency. This does not contradict Corollary~\ref{corol:strong-implies-onestep} because the notion of consistency in~\cite{TAWJ} is different from Definition~\ref{def:consistency} and seems more appropriate for the dynamics of these networks. Essentially,  switches of the Boolean states of several variables are treated as simultaneous if they occur within a small time-windows that mark the boundaries between so-called \emph{episodes.}  This opens a promising avenue of future research on generalizations of this notion of consistency as well as possible generalizations of Theorem~\ref{thm:mainth} in this direction.

Gehrmann and Drossel~\cite{Gehrmann} obtained results about consistency between ODE models for two small gene regulatory networks and their Boolean approximations.  In contrast to Glass networks and~\cite{TAWJ}, in these models the right-hand sides of the ODEs are Lipschitz-continuous, as they are in the models studied here.

The  systems of~\cite{Gehrmann, TAWJ} and our systems~$\bD_2(f, \vg)$ can all be conceptualized as
\emph{networks of interacting agents,} where the ODE variables are grouped into pairwise disjoint sets~$X_i$ that constitute the \emph{$i$-th agent} and only one Boolean variable is assigned to each agent. Roughly speaking, in~\cite{TAWJ} an agent corresponds to a neuron whose state is characterized by a cross-membrane voltage and a so-called gating variable, in~\cite{Gehrmann} an agent corresponds to a gene whose state is characterized by its mRNA concentration and the concentration of the corresponding gene product, and in~$\bD_2(f, \vg)$ we have
$X_i = \{x_i, x_{i+n}\}$.  Such a partition into agents appears natural in many models of real-world systems and may favor consistency between an ODE system and its Boolean approximation.  For example, in our systems~$\bD_2(f, \vg)$ self-regulation does not cause the kind of problems for consistency as it may cause in systems~$\bD_1(f, \vg)$  or Glass networks, where all variables are discretized.   Identifying and perhaps classifying general mechanisms in such networks of interacting agents that favor or guarantee some form of consistency is another promising avenue for further research.

In particular, we are planning to generalize Theorem~\ref{thm:mainth} to much larger classes of networks of interacting agents.
The mechanism that gives Theorem~\ref{thm:mainth} of course is sufficient separation of timescales in addition to suitable bifurcations in the internal dynamics of the signature variables, and the theorem should generalize in this form. The consistency result in~\cite{TAWJ} also relies on separation of time scales, with voltages changing much faster than gating variables, but an additional synchronization mechanism is built into the network architecture. In contrast to Theorem~\ref{thm:mainth}, very large separation of time scales destroys the consistency for the first system studied in~\cite{Gehrmann}, despite the fact that its Boolean approximation is one-stepping.  Thus the mechanism for consistency at other time scales that is reported in~\cite{Gehrmann} must be different.

There are some similarities between the proof of our Theorem~\ref{thm:mainth} and the results in~\cite{Ironi} about certain models of gene regulatory networks with steep sigmoid functions.  In particular, we want to note that in our systems
$\bD_2(f, \vg)$ each signaling variable regulates exactly one signature variable has a parallel in Assumption~$\mathcal{A}$ of~\cite{Ironi} that every transcription factor regulates exactly one gene at each of its thresholds.  It appears, however, that proving analytical results about gene regulation with steep continuous sigmoid function is in general more challenging than proving analogous results for the classes described here (see~\cite{Ironi} and references therein for description of some difficulties in the former classes).

Many Boolean systems that have been studied in the literature are one-stepping.  The \emph{one-input systems} considered in~\cite{GlassKauffman} are a special case, we mentioned already~\cite{Gehrmann}, and the Boolean approximation of the models considered in~\cite{Polynikis} are another example. The latter paper compares various approximations, not only Boolean ones, to a complete nonlinear model of gene regulation and reports that substantial discrepancies can arise between the predictions of various models for even such a simple system.  This is another piece of evidence that our models are better behaved than would be expected of models of gene regulation.

\subsection{Chaos}\label{subsec:chaos}

Chaos in differentiable flows has been widely studied.  A notion of chaos in classes of Boolean systems has been proposed by S.~Kauffman and has been fairly extensively studied (see~\cite{Aldana, origins} for  reviews).  It is a natural question whether there is any connection between these two notions of chaos in the sense that differentiable flows that can be translated into a chaotic  Boolean system must be or will tend to be chaotic flows, or \emph{vice versa.} In other words, it is natural to ask whether the two notions of chaos are saying something equivalent about the underlying natural system or whether they are just two different mathematical properties that happen to share a name.

Chaos in individual Boolean systems can be defined in terms of the \emph{slope of the Derrida curve at the origin,} that is, the average Hamming distance after one updating step for randomly chosen initial states with Hamming distance~1. If this slope exceeds~1, the system can be considered chaotic. It is easy to see that there cannot be a straightforward correspondence between chaos in flows and chaos in this sense in their Boolean approximations.
For example, the two-dimensional system with~$f$ given by

$$f(00) = 00, \quad f(01) = 11, \quad f(10) = 01, \quad f(11) = 10$$
is chaotic in this sense, but the corresponding two-dimensional flow
$\bD_1(f, \vg)$ cannot exhibit chaos.

It seems quite likely though that there are some more subtle connections between the two notions of chaos or between chaos and consistency. Exploring such connections leads to interesting and quite challenging questions.  For example, it would be interesting to know  under what conditions the sets $U^s$ of Theorem~\ref{thm:mainth}(ii) must intersect a unique (attracting) periodic orbit, which would preclude chaos in this region of the state space.  This remains an open question even for the one-stepping trajectories of Example~\ref{ex:fixed-vs-onestep}(ii) and the conversion schemes that we presented in Subsection~\ref{subsec:convert-schemes}. Our software has a module for estimating Lyapunov exponents and allows for numerical explorations of such problems.

\subsection{Conclusions}

We constructed, for every Boolean system with updating function~$f$, two classes of ODE systems $\bD_1(f, \vg, \cQ), \bD_2(f, \vg, \cQ)$ for which the Boolean system
$\bB = (2^n, f)$ is a natural approximation. Theorem~\ref{thm:mainth} shows that sufficiently large separation of time scales guarantees consistency between
$\bD_2(f, \vg, \cQ)$ and $\bB = (2^n, f)$ on a large region of the state space of the former, as long as~$f$ is monotone-stepping. Classes of ODE systems that were previously investigated in the literature for consistency with their natural Boolean approximations either have discontinuities in their right-hand sides or tend to not satisfy the counterpart of Theorem~\ref{thm:mainth} in its most general form. In our opinion, this makes the models presented here a promising class of toy models for elucidating general mechanisms that favor or entail consistency between an ODE system and its Boolean approximation. It seems likely that the insights from studying such mechanisms in the context of our models would, at least to some extent, carry over to classes of realistic but analytically less tractable models of natural, in particular biological, phenomena.

\section*{Acknowledgement}

We thank Bismark Oduro and Hanyan Zhu for valuable contributions to the preliminary stages of the project that resulted in this paper.

\end{document}